\documentstyle[11pt]{article}

\author{Charles H. Walter\thanks{Supported in part by NSA research grant
MDA904-92-H-3009.}\\
Laboratoire de Math\'ematiques\\
Universit\'e de Nice\\
F-06108 Nice Cedex 02 FRANCE}
\title{Irreducibility of Moduli Spaces of Vector Bundles on Birationally Ruled
Surfaces
}
\date{ }

\input amssym.def
\newsymbol\ncong 231D
\newsymbol\boxtimes 1202
\newsymbol\twoheadrightarrow 1310

\newtheorem{theorem}{Theorem}
\newtheorem{lemma}[theorem]{Lemma}
\newtheorem{proposition}[theorem]{Proposition}

\def\qed{\hfill $\Box$}
\long\def\TeXButton#1#2{#2}
\def\proof{\paragraph{Proof. }}

\catcode`\@=11
\def\@begintheorem#1#2{\sl \trivlist
   \item[\hskip \labelsep{\bf #1\ #2\thmcounterend}]}
\def\@opargbegintheorem#1#2#3{\sl \trivlist
      \item[\hskip \labelsep{\bf #1\ #2\ (#3)\thmcounterend}]}
\def\thmcounterend{.}

\def\section{\@startsection{section}{1}{\z@}{-3.25ex plus
 -1ex minus -.2ex}{1.5ex plus .2ex}{\large\bf}}
\def\subsection{\@startsection
 {subsection}{2}{\z@}{3.25ex plus 1ex minus .2ex}{-0.5em}{\normalsize\sl}}
\def\subsubsection{\@startsection
 {subsubsection}{3}{\z@}{3.25ex plus 1ex minus .2ex}{-0.5em}{\normalsize\sl}}
\def\paragraph{\@startsection
 {paragraph}{3}{\z@}{2ex plus 0.6ex minus .2ex}{-0.5em}{\normalsize\sl}}
\def\subparagraph{\@startsection
 {subparagraph}{3}{\parindent}{2ex plus 0.6ex minus
.2ex}{-1pt}{\normalsize\sl}}

\catcode`\@=12

\begin{document}

\maketitle
\begin{abstract}
\noindent Let $S$ be a birationally ruled surface. We show that the moduli
schemes $M_S(r,c_1,c_2)$ of semistable sheaves on $S$ of rank $r$ and Chern
classes $c_1$ and $c_2$ are irreducible for all $(r,c_1,c_2)$ provided the
polarization of $S$ used satisfies a simple numerical condition. This is
accomplished by proving that the stacks of prioritary sheaves on $S$ of
fixed rank and Chern classes are smooth and irreducible.\bigskip\
\end{abstract}

One important recent result in the theory of vector bundles on algebraic
surfaces is the theorem of Gieseker and Li that for any smooth projective
surface $S$ and any ample divisor $H$ on $S$, the moduli scheme $%
M_{S,H}(2,c_1,c_2)$ of S-equivalence classes of $H$-semistable torsion-free
sheaves of rank $2$, determinant $c_1\in {\rm Pic}(S)$, and second Chern
class $c_2$ is irreducible if $c_2\gg 0$. If $S$ is a surface of general
type, the condition $c_2\gg 0$ is necessary because of an example of
Gieseker with small $c_2$ where the moduli space is reducible. In contrast
it has been known for quite some time that the moduli schemes $M_{{\bf P}%
^2,H}(r,c_1,c_2)$ is irreducible for all $(r,c_1,c_2)$ for which there exist
semistable sheaves on the projective plane, and the same result is also
known for ${\bf P}^1\times {\bf P}^1$.

In this paper we extend this strong irreducibility result from ${\bf P}^2$
and ${\bf P}^1\times {\bf P}^1$ to all smooth projective surfaces of
negative Kodaira dimension. To simplify our exposition, we will omit ${\bf P}%
^2$ although it can be handled by the same method. (Indeed our method is
based on a method of Ellingsrud and Str\o mme which was developed for ${\bf P%
}^2$.) So our surface $S$ possesses\ a morphism $\pi {:}~S\rightarrow C$
onto a smooth curve with connected fibers and with general fiber isomorphic
to ${\bf P}^1$. We fix such a $\pi $. (Such a $\pi $ is unique if $%
q(S)=g(C)>0$ or if $S={\bf P}({\cal O}_{{\bf P}^1}\oplus {\cal O}_{{\bf P}%
^1}(e))$ with $e>0$, but there can be many, even infinitely many, possible $%
\pi $ for certain rational surfaces.) For $p\in C$, let $f_p=\pi ^{-1}(p)$.
These $f_p$ are all numerically equivalent, and we write $f\in {\rm NS}(X)$
for the numerical class of these $f_p$. We prove

\begin{theorem}
\label{semistable}Let $\pi {:}~S\rightarrow C$ be a birationally ruled
surface and $f\in {\rm NS}(S)$ the numerical class of a fiber of $\pi $. Let
$H$ be an ample divisor on $S$ such that $H\cdot (K_S+f)<0$. Suppose $r\geq 2
$, $c_1\in {\rm NS}(S)$, and $c_2\in {\bf Z}$ are given. If the moduli
scheme $M_{S,H}(r,c_1,c_2)$ of S-equivalence classes of $H$-semistable
torsion-free sheaves of rank $r$ and Chern classes $c_1$ and $c_2$ is
non-empty, then it is irreducible and normal. In addition, the open
subscheme $M_{S,H}^s(r,c_1,c_2)$ parametrizing stable sheaves is smooth.
\end{theorem}

Our methods also show that the general $H$-semistable torsion-free sheaf in
any of the $M_{S,H}(r,c_1,c_2)$ is locally free and that there is a
dominant, generically finite map from an open subscheme of ${\rm Jac}%
(C)\times {\rm Jac}(C)\times {\bf P}^m$ to $M_{S,H}^s(r,c_1,c_2)$ with $%
m=2rc_2-(r-1)c_1+(r^2-2)g-r^2+1$ where $g$ is the genus of $C$. So if $S$ is
a rational surface, then $M_{S,H}^s(r,c_1,c_2)$ is unirational.

Ample divisors $H$ satisfying the hypothesis $H\cdot (K_S+f)<0$ exist on any
birationally ruled surface because of Lemma \ref{H} below. On certain
surfaces there may exist a divisor $D$ of degree $1$ on $C$ such that the
divisor class $-K_S-\pi ^{*}(D)$ is effective. In that case all ample
divisors $H$ satisfy $H\cdot (K_S+f)<0$, and the theorem holds for all
possible polarizations. Examples of such surfaces include Del Pezzo
surfaces, rational ruled surfaces, and ruled surfaces of the form ${\bf P}(%
{\rm O}_C\oplus {\rm O}_C(D_0))$ with $D_0$ a divisor of degree at least $%
2g-1$ on the smooth projective curve $C$ of genus $g$.

Our method of proof begins by adapting a definition from \cite{HL}. If $\pi {%
:}~S\rightarrow C$ is a birationally ruled surface, then we will say that a
coherent sheaf ${\cal E}$ on $S$ is {\em prioritary} (with respect to $\pi $%
) if it is torsion-free and if ${\rm Ext}^2({\cal E},{\cal E}(-f_p))=0\,$
for all $p\in C$. By the semicontinuity theorem the prioritary sheaves in
any locally noetherian flat family of coherent sheaves on $S$ form an open
subfamily. Hence the prioritary sheaves on $S$ are parametrized by an open
substack of the stack of coherent sheaves on $S$.

For a given $r\geq 1$, $c_1\in {\rm NS}(S)$, and $c_2\in {\bf Z}$, we will
write ${\rm Coh}_S(r,c_1,c_2)$ for the stack of coherent sheaves of rank $r$
and Chern classes $c_1$ and $c_2$ (modulo numerical equivalence), and ${\rm %
TF}_S(r,c_1,c_2)$ and ${\rm \Pr ior}_S(r,c_1,c_2)$ for the open substacks
of, respectively, torsion-free and prioritary sheaves. We will derive
Theorem \ref{semistable} from:

\begin{proposition}
\label{stack}Let $\pi {:}~S\rightarrow C$ be a birationally ruled surface.
Suppose $r\geq 2$, $c_1\in {\rm NS}(S)$, and $c_2\in {\bf Z}$ are given.
Then the stack ${\rm \Pr ior}_S(r,c_1,c_2)$ of prioritary sheaves on $S$ of
rank $r$ and Chern classes $c_1$ and $c_2$ is smooth and irreducible.
\end{proposition}

The Proposition is proven in two steps. First we prove it for geometrically
ruled surfaces. Then we prove that for $S\rightarrow S_1$ the blowup of a
point, if the proposition holds for $S_1$ then it holds for $S$. Our method
is based on the version of the method of Ellingsrud and Str\o mme (\cite{E},
\cite{ES}, \cite{HS}) as presented in \cite{LP} \S 9.

The reader unfamiliar with algebraic stacks may wish to consult \cite{LMB}.
We use stacks because in that context there exist natural universal families
of coherent (or torsion-free or prioritary) sheaves. Alternative universal
families which stay within the category of schemes would be certain standard
open subschemes of Quot schemes. But these depend on the choice of a
polarization ${\cal O}_S(1)$, of the Hilbert polynomial $P$, of twists $m\gg
0$ and of a vector space $H_m$ of dimension $P(m)$. One then deals with the
scheme ${\rm Quot}_{S,{\cal O}_S(1)}^0(P,m)$ parametrizing all quotients $%
\gamma {:}~H_m\otimes {\cal O}_S(-m)\TeXButton{-->>}{\twoheadrightarrow}%
{\cal F}$ such that $H^i({\cal F}(m))=0$ for all $i\geq 1$ and such that the
induced map $H_m\rightarrow H^0({\cal F}(m))$ is an isomorphism. One would
prefer not to work in a context where one constantly has to refer to all
these choices, particularly since no single set of choices will work when
one deals with unlimited families. But strictly speaking, these Quot schemes
are not that far from our point of view because the verification in \cite
{LMB} (4.14.2) that the coherent sheaves on $S$ are parametrized by an
algebraic stack ${\rm Coh}_S$ is done essentially by gluing together all the
${\rm Quot}_{S,{\cal O}_S(1)}^0(P,m)$ in the smooth Grothendieck topology.

This paper was written in the context of the group on vector bundles on
surfaces of Europroj and as a direct result of the Catania congress where
the problem was mentioned by J.\ Le Potier. The author would also like to
thank A.\ Hirschowitz for some pertinent comments.

\section{Proof of the Theorem}

We begin with two lemmas about coherent sheaves on ${\bf P}^1$ and the
restriction of torsion-free sheaves on surfaces to curves in the surface.
These lemmas are well known although they have usually been stated in terms
of complete families or of versal deformation spaces instead of stacks. We
state them without proof.

\begin{lemma}
\label{P1}Let $r\geq 2$ and $0\leq d<r$ be integers. Let ${\rm Coh}_{{\bf P}%
^1}(r,-d)$ be the stack of coherent sheaves of rank $r$ and degree $-d$ on $%
{\bf P}^1$.

(i)\quad If $d>0$, then sheaves not isomorphic to ${\cal O}_{{\bf P}%
^1}^{r-d}\oplus {\cal O}_{{\bf P}^1}(-1)^d$ form a closed substack of ${\rm %
Coh}_{{\bf P}^1}(r,-d)$ of codimension at least $2$.

(ii)\quad If $d=0$, then sheaves not isomorphic to ${\cal O}_{{\bf P}^1}^r$
form a closed substack of ${\rm Coh}_{{\bf P}^1}(r,0)$ of codimension $1$.
Sheaves isomorphic neither to ${\cal O}_{{\bf P}^1}^r$ nor to ${\cal O}_{%
{\bf P}^1}(1)\oplus {\cal O}_{{\bf P}^1}^{r-2}\oplus {\cal O}_{{\bf P}^1}(-1)
$ form a closed substack of ${\rm Coh}_{{\bf P}^1}(r,0)$ of codimension at
least $2$.
\end{lemma}

\begin{lemma}
\label{restriction}Let $D$ be an effective Cartier divisor on a projective
surface $S$. If ${\cal E}$ is a torsion-free sheaf on $S$ such that ${\rm Ext%
}^2({\cal E},{\cal E}(-D))=0$, then the restriction map ${\rm TF}%
_S(r,c_1,c_2)\rightarrow {\rm Coh}_D(r,c_1\cdot D)$ is smooth (and therefore
open) in a neighborhood of $[{\cal E}]$.
\end{lemma}

We also need two lemmas for reduction steps in the proof of Proposition \ref
{stack}.

\begin{lemma}
\label{Beil}Let $\pi {:}~S\rightarrow C$ be a geometrically ruled surface
with a section $\sigma \subset S$. If ${\cal E}$ is a coherent sheaf on $S$
such that $\pi _{*}({\cal E}(-\sigma ))=R^1\pi _{*}({\cal E})=0$, then there
is an exact sequence%
$$
0\rightarrow \pi ^{*}(\pi _{*}({\cal E}))\rightarrow {\cal E}\rightarrow \pi
^{*}(R^1\pi _{*}({\cal E}(-\sigma )))\otimes \Omega _{S/C}(\sigma
)\rightarrow 0.
$$
\end{lemma}

\TeXButton{Proof}{\proof}This is a special case of a relative version of
Beilinson's spectral sequence, but for lack of a precise reference we give
the proof in full. Let $Y:=S\times _CS$. Then the diagonal $\Delta $ of $Y$
has Beilinson's resolution (cf.\ \cite{B})%
$$
0\rightarrow \Omega _{S/C}(\sigma )\TeXButton{boxtimes}{\boxtimes}{\cal O}%
_S(-\sigma )\rightarrow {\cal O}_Y\rightarrow {\cal O}_\Delta \rightarrow 0.
$$
Applying $R^i{\rm pr}_{1*}(-\otimes {\rm pr}_2^{*}({\cal E}))$ to this exact
sequence gives a long exact sequence which is equivalent to the one asserted
by the lemma because one always has $R^i{\rm pr}_{1*}({\cal F}
\TeXButton{boxtimes}{\boxtimes}{\cal G})\cong {\cal F}\otimes \pi
^{*}(R^i\pi _{*}({\cal G}))$ if ${\cal F}$ is locally free and ${\cal G}$
coherent on $S$ because of the projection formula and \cite{H}, Chapter III,
Proposition 9.3. \TeXButton{qed}{\qed \medskip}

\begin{lemma}
\label{blowup}Let $S_1$ be a smooth surface $\alpha {:}~S\rightarrow S_1$
the blowup of a point $x$ of $S_1$. Let $E$ be the exceptional divisor in $S$%
. Suppose that ${\cal E}$ is a coherent sheaf of rank $r$ on $S$ such that $%
{\cal E}{\mid }_E\cong {\cal O}_E^{r-d}\oplus {\cal O}_E(-1)^d$ for some $d$%
. Then $\alpha _{*}({\cal E})$ is locally free in a neighborhood of $x$, and
there are exact sequences%
\begin{eqnarray*}
& 0\rightarrow \alpha ^{*}(\alpha _{*}({\cal E}))\rightarrow {\cal E}%
\rightarrow {\cal O}_E(-1)^d\rightarrow 0, &  \\
& 0\rightarrow {\cal E}(-E)\rightarrow \alpha ^{*}(\alpha _{*}({\cal E}%
))\rightarrow {\cal O}_E^{r-d}\rightarrow 0. &
\end{eqnarray*}
Moreover, for any divisor $D$ on $S_1$ we have ${\rm Ext}^2({\cal E},{\cal E}%
(\alpha ^{*}(D)))\cong {\rm Ext}^2(\alpha _{*}({\cal E}),\alpha {\cal _{*}(E}%
)(D))$.
\end{lemma}

\TeXButton{Proof}{\proof}Let ${\cal F}$ be the kernel of the composition $%
{\cal E}\TeXButton{-->>}{\twoheadrightarrow}{\cal E}{\mid }_E\TeXButton{-->>}
{\twoheadrightarrow}{\cal O}_E(-1)^d$. By general properties of elementary
transforms the exact sequence $0\rightarrow {\cal O}_E^{r-d}\rightarrow
{\cal E}{\mid }_E\rightarrow {\cal O}_E(-1)^d\rightarrow 0$ transforms into $%
0\rightarrow {\cal O}_E^d\rightarrow {\cal F}{\mid }_E\rightarrow {\cal O}%
_E^{r-d}\rightarrow 0$. So ${\cal F}$ is trivial along $E$, and ${\cal F}%
\cong \alpha ^{*}(\alpha _{*}({\cal F}))$. Applying $\alpha _{*}$ to the
exact sequence $0\rightarrow {\cal F}\rightarrow {\cal E}\rightarrow {\cal O}%
_E(-1)^d\rightarrow 0$, we see that $\alpha _{*}({\cal F})\cong \alpha _{*}(%
{\cal E})$. Hence ${\cal F}\cong \alpha ^{*}(\alpha _{*}({\cal E}))$. The
exact sequences asserted by the lemma are now the standard exact sequences
of an elementary transform.

By adjunction and the formula $K_S=\alpha ^{*}(K_{S_1})+E$ we see that%
$$
{\rm Hom}(\alpha _{*}({\cal E}),\alpha {\cal _{*}(E})(-D+K_{S_1}))\cong {\rm %
Hom}({\cal F},{\cal E}(-\alpha ^{*}(D)+K_S-E))
$$
or by Serre duality that ${\rm Ext}^2(\alpha _{*}({\cal E}),\alpha {\cal %
_{*}(E})(D))\cong {\rm Ext}^2({\cal E},{\cal F}(\alpha ^{*}(D)+E))$. If we
now apply the functor ${\rm Ext}^2({\cal E}(-\alpha ^{*}(D)),-)$ to the
exact sequence $0\rightarrow {\cal E}\rightarrow {\cal F}(E)\rightarrow
{\cal O}_E(-1)^{r-d}\rightarrow 0$ and note that
$$
{\rm Ext}^i({\cal E}(-\alpha ^{*}(D)),{\cal O}_E(-1))\cong H^i(E,({\cal E}{%
\mid }_E)^{\lor }(-1))=0
$$
for $i=1,2$, then we see that ${\rm Ext}^2({\cal E},{\cal F}(\alpha
^{*}(D)+E))\cong {\rm Ext}^2({\cal E},{\cal E}(\alpha ^{*}(D)))$. This
completes the proof of the lemma. \TeXButton{qed}{\qed \medskip}\

We now begin the proof of Proposition \ref{stack}. The smoothness of ${\rm %
\Pr ior}_S(r,c_1,c_2)$ follows from ${\rm Ext}^2({\cal E},{\cal E})=0$
because this is the obstruction space for deformations of ${\cal E}$. So we
concentrate on irreducibility. We begin with the special case of
geometrically ruled surfaces.

\paragraph{Proof of Proposition \ref{stack} when $\pi {:}~S\rightarrow C$ is
a geometrically ruled surface.}

We follow the method of Ellingsrud and Str\o mme as presented in \cite{LP}
\S 9.

We fix a section $\sigma \subset $$S$. Replacing ${\cal E}$ by an
appropriate twist ${\cal E}(n\sigma )$ if necessary we may assume that $%
d:=-c_1\cdot f$ satisfies $0\leq d<r$. The proof now divides briefly into
two cases $d>0$ and $d=0$.

If $d>$$0$, then by Lemmas \ref{P1} and \ref{restriction}, those ${\cal E}$
such that ${\cal E}{\mid }_{f_p}\TeXButton{ncong}{\ncong}{\cal O}%
_{f_p}^{r-d}\oplus {\cal O}_{f_p}(-1)^d$ for some $p\in C$ are parametrized
by a closed substack of ${\rm \Pr ior}_S(r,c_1,c_2)$ of codimension at least
$1$. So we may restrict ourselves to the dense open substack ${\rm \Pr ior}%
^0 $ where ${\cal E}{\mid }_{f_p}\cong {\cal O}_{f_p}^{r-d}\oplus {\cal O}%
_{f_p}(-1)^d$ for all $p\in C$.

If $d=0$, then by an analogous argument, we may restrict ourselves to a
dense open substack ${\rm \Pr ior}^0$ where ${\cal E}{\mid }_{f_p}\cong
{\cal O}_{f_p}^r$ for all $p\in C$ except for a finite number of $p$ where $%
{\cal E}{\mid }_{f_p}\cong {\cal O}_{f_p}(1)\oplus {\cal O}%
_{f_p}^{r-2}\oplus {\cal O}_{f_p}(-1)$.

In either case, we set ${\cal K}:=\pi _{*}({\cal E)}$. Since $R^1\pi _{*}(%
{\cal E})=0$, the Leray-Serre spectral sequence implies that $\chi ({\cal E}%
)=\chi ({\cal K})$. We may now calculate by Riemann-Roch that ${\cal K}$ is
a vector bundle on $C$ of rank $r-d$ and degree $k:=\chi ({\cal E}%
)+(r-d)(g-1)$ where $g$ is the genus of $C$.

Let ${\cal L}=R^1\pi _{*}({\cal E}(-\sigma ))$. Then ${\cal L}$ is a sheaf
on $C$ of rank $d$ and degree $l:=-\chi ({\cal E}(-\sigma ))+d(g-1)=-\chi (%
{\cal E})+(c_1\cdot \sigma )-(r-d)(g-1)$. The sheaf ${\cal L}$ is locally
free if $d>0$.

By Lemma \ref{Beil} there is an exact sequence%
$$
0\rightarrow \pi ^{*}({\cal K})\rightarrow {\cal E}\rightarrow \pi ^{*}(%
{\cal L})\otimes \Omega _{S/C}(\sigma )\rightarrow 0.
$$
Now using the notations ${\rm Ext}_{+}$ and ${\rm Ext}_{-}$ of \cite{DLP},
p.\ 200, we have%
$$
{\rm Ext}_{+}^i({\cal E},{\cal E})=H^i(\pi ^{*}({\cal K}^{\vee }\otimes
{\cal L})\otimes \Omega _{S/C}(\sigma ))=0
$$
for all $i$. Hence ${\rm Ext}^i({\cal E},{\cal E})\cong {\rm Ext}_{-}^i(%
{\cal E},{\cal E})$ for all $i$. Thus the infinitesimal deformations of $%
{\cal E}$ are the same as the infinitesimal deformations of the filtered
sheaf  $0\subset \pi ^{*}({\cal K})\subset {\cal E}$. Furthermore, since $%
{\rm Ext}^2(\pi ^{*}({\cal L})\otimes \Omega _{S/C}(\sigma ),\pi ^{*}({\cal K%
}))=0$, we have a surjection%
$$
{\rm Ext}_{-}^1({\cal E},{\cal E})\rightarrow {\rm Ext}_{{\cal O}_C}^1({\cal %
K},{\cal K})\oplus {\rm Ext}_{{\cal O}_C}^1({\cal L},{\cal L})\rightarrow 0.
$$
Hence a general infinitesimal deformation of ${\cal E}$ induces general
infinitesimal deformations of ${\cal K}$ and ${\cal L}$. Since none of these
deformations are obstructed, the morphism ${\rm \Pr ior}^0\rightarrow {\rm %
Coh}_C(r-d,k)\times {\rm Coh}_C(d,l)$ defined by $[{\cal E}]\mapsto ([{\cal K%
}],[{\cal L}])$ is dominant. The fibers of this morphism are irreducible
since they are stack quotients of an open subscheme of the affine space $%
{\rm Ext}^1(\pi ^{*}({\cal L})\otimes \Omega _{S/C}(\sigma ),\pi ^{*}({\cal K%
}))$.The target of the morphism is irreducible since stacks of coherent
sheaves of a fixed rank and degree on a smooth connected curve $C$ are
irreducible. So ${\rm \Pr ior}^0$ and hence also ${\rm \Pr ior}_S(r,c_1,c_2)$
are irreducible. This completes of the proof of Proposition \ref{stack} when
$\pi {:}~S\rightarrow C$ is a geometrically ruled surface. \TeXButton{qed}
{\qed \medskip}

\paragraph{Proof of Proposition \ref{stack} in general.}

We go by induction on the Picard number $\rho (S):={\rm rk}_{{\bf Z}}({\rm NS%
}(S))$. The initial value is $\rho (S)=2$ which is the case of geometrically
ruled surfaces which we just proved. So we may assume that $\rho (S)\geq 3$.
Then $\pi {:}~S\rightarrow C$ is birationally but not geometrically ruled.
So some fiber of $\pi $ contains an irreducible component $E\cong {\bf P}^1$
such that $E^2=-1$. We let $\alpha {:}~S\rightarrow S_1$ be the contraction
of $E$, and we let $\beta {:}~S_1\rightarrow C$ be the morphism such that $%
\pi $ factors as $\pi =\beta \alpha $. Since $\rho (S_1)=\rho (S)-1$, we may
assume by induction that Proposition \ref{stack} holds for $\beta {:}%
{}~S_1\rightarrow C$.

Let $d=-c_1\cdot E$. By replacing ${\cal E}$ with an appropriate twist $%
{\cal E}(nE)$ we may assume that $0\leq d<r$. Because $f_{\pi (E)}-E$ is
effective, the condition ${\rm Ext}^2({\cal E},{\cal E}(-f_{\pi (E)}))=0$
implies ${\rm Ext}^2({\cal E},{\cal E}(-E))=0$. So by Lemmas \ref{P1} and
\ref{restriction} the substack ${\rm \Pr ior}^1\subset {\rm \Pr ior}%
_S(r,c_1,c_2)$ parametrizing prioritary sheaves ${\cal E}$ such that ${\cal E%
}{\mid }_E\cong {\cal O}_E^{r-d}\oplus {\cal O}_E(-1)^d$ is open and dense.
By Lemma \ref{blowup}, the application $[{\cal E}]\rightarrow [\alpha _{*}(%
{\cal E})]$ defines a morphism ${\rm \Pr ior}^1\rightarrow {\rm \Pr ior}%
_{S_1}(r,\alpha _{*}(c_1),c_2+\frac 12d(d-1))$. Moreover, the morphism
realizes ${\rm \Pr ior}^1$ as a $d(r-d)$-dimensional Grassmannian bundle
over the dense open substack of ${\rm \Pr ior}_{S_1}(r,\alpha
_{*}(c_1),c_2+\frac 12d(d-1))$ of sheaves which are locally free at the
center of the blowup. Since ${\rm \Pr ior}_{S_1}(r,\alpha
_{*}(c_1),c_2+\frac 12d(d-1))$ is irreducible by the inductive hypothesis,
we see that ${\rm \Pr ior}^1$ and hence also ${\rm \Pr ior}_S(r,c_1,c_2)$
are irreducible. This completes the proof of Proposition \ref{stack}.
\TeXButton{qed}{\qed \medskip}

\paragraph{Proof of Theorem \ref{semistable}.}

First we show that if $H\cdot (K_S+f)<0$, then any $H$-semistable sheaf $%
{\cal E}$ is prioritary. But if ${\cal E}$ is $H$-semistable, then any
nonzero torsion-free quotient ${\cal Q}$ of ${\cal E}$ would have $H$-slope
satisfying $\mu _H({\cal Q})\geq \mu _H({\cal E})$, while any nonzero
subsheaf ${\cal S}$ of ${\cal E}$ would have $H$-slope satisfying $\mu _H(%
{\cal S})\leq \mu _H({\cal E})$. So if ${\cal E}$ were not prioritary, there
would exist a $p\in C$ such that ${\rm Ext}^2({\cal E},{\cal E}(-f_p))\neq 0$%
. There would then exist a nonzero $\phi \in {\rm Hom}({\cal E},{\cal E}%
(K_S+f_p))\cong {\rm Ext}^2({\cal E},{\cal E}(-f_p))^{*}$. The image of $%
\phi $ would then satisfy%
$$
\mu _H({\cal E})\leq \mu _H({\rm im}(\phi ))\leq \mu _H({\cal E}%
(K_S+f_p))=\mu _H({\cal E})+H\cdot (K_S+f),
$$
contradicting $H\cdot (K_S+f)<0$.

Thus the semistable sheaves on $S$ of rank $r$ and Chern classes $c_1$ and $%
c_2$ are parametrized by an open substack $H{\rm -SS}\subset {\rm \Pr ior}%
_S(r,c_1,c_2)$. This last stack is smooth and irreducible by Proposition \ref
{stack}. So if there exist $H$-semistable sheaves on $S$ with that rank and
those Chern classes, then $H{\rm -SS}$ is a smooth and irreducible stack.

We now show that this implies that the moduli scheme $M_{S,H}(r,c_1,c_2)$ is
normal and that the open subscheme $M_{S,H}^s(r,c_1,c_2)$ parametrizing $H$%
-stable sheaves is smooth. We write ${\cal O}_S(1):={\cal O}_S(H)$. Since $H%
{\rm -SS}$ is a limited family, there exists an integer $m\gg 0$ such that
all $H$-semistable sheaves ${\cal E}$ in $H{\rm -SS}$ satisfy $H^i({\cal E}%
(m))=0$ for $i\geq 1$ and have ${\cal E}(m)$ generated by global sections.
Let $H_m$ be a vector space of dimension $h^0{\cal (E}(m))$ and let $Q^{ss}:=%
{\rm Quot}_{S,H}^{ss}(m;r,c_1,c_2)$ denote the Hilbert-Grothendieck scheme
parametrizing all quotients $\gamma {:}~H_m\otimes {\cal O}_S(-m)
\TeXButton{-->>}{\twoheadrightarrow}{\cal F}$ such that ${\cal F}$ is $H$%
-semistable of rank $r$ and Chern classes $c_1$and $c_2$ with $H^i({\cal F}%
(m))=0$ for all $i\geq 1$ and such that the induced map $H_m\rightarrow H^0(%
{\cal F}(m))$ is an isomorphism. Then according to the construction of \cite
{LMB} (4.14.2), $H{\rm -SS}$ may be identified with the quotient stack$\
\left[ Q^{ss}/{\rm GL}(H_m)\right] $. Hence $Q^{ss}$ is a smooth and
irreducible scheme. But $M_{S,H}(r,c_1,c_2)$ is the GIT quotient scheme $%
Q^{ss}//({\rm SL}(H_m),{\cal L}_N)$ where ${\cal L}_N$ is Simpson's
polarization of $Q^{ss}$ defined by ${\cal L}_N:=\det ({\rm pr}_{1*}({\cal U}%
\otimes {\rm pr}_2^{*}({\cal O}_S(N))))$ where ${\cal U}$ is the universal
sheaf on $Q^{ss}\times S$, the ${\rm pr}_i$ are the two projections, and $%
N\gg m$ (cf.\ \cite{S} \S 1). Thus $M_{S,H}(r,c_1,c_2)$ is the GIT quotient
of a smooth and irreducible variety. But such quotients, when nonempty, are
always normal and irreducible varieties, and the points of the quotient
corresponding to stable points are smooth. \TeXButton{qed}{\qed \medskip}

\section{Existence of Ample Divisors}

We show that on any birationally ruled surface there exists an ample divisor
satisfying the hypothesis of Theorem \ref{semistable}.

\begin{lemma}
\label{H}Let $S$ be a birationally ruled surface, $\pi {:}~S\rightarrow C$ a
birational ruling, and $f\in {\rm NS}(X)$ the class of a fiber of $\pi $.
Then there exists an ample divisor $H$ on $S$ such that $H\cdot (K_S+f)<0$.
\end{lemma}

\TeXButton{Proof}{\proof}First we consider the case of $S$ a geometrically
ruled surface. Let $\sigma $ be a section of minimal self-intersection $-e$.
According to \cite{H} Chapter V, Corollary 2.11, Propositions 2.20 and 2.21
and Exercise 2.14, we see that $K_S\equiv -2\sigma +(2g-2-e)f$ and that an $%
H=\sigma +bf$ is ample if $b$ is sufficiently large. Since $H\cdot
(K_S+f)=2g-1+e-2b$, by picking $b$ large enough we get an ample $H$ such
that $H\cdot (K_S+f)<0$.

Now we consider the case where $\pi {:}~S\rightarrow C$ is birationally but
not geometrically ruled. Then some fiber of $\pi $ contains an exceptional
divisor of the first kind $E$. Let $\alpha :S\rightarrow S_1$ be the
contraction of $E$. By induction on the Picard number we may assume there
exists an ample divisor $H_1$ on $S_1$ such that $H_1\cdot (K_{S_1}+f)<0$.
{}From \cite{H} Chapter V, Proposition 3.3 and Exercise 3.3, we see that $%
K_S=\alpha ^{*}(K_{S_1})+E$ and that $H:=2\alpha ^{*}(H_1)-E$ is ample. Then
$H\cdot (K_S+f)=2(H_1\cdot (K_{S_1}+f))-E^2<0$. \TeXButton{qed}
{\qed \medskip}

\end{document}